\documentclass[12pt]{article}

\let\oldmarginpar\marginpar
\renewcommand\marginpar[1]{\-\oldmarginpar[\raggedleft\footnotesize #1]%
{\raggedright\footnotesize #1}}

\usepackage{mathptmx}
\usepackage{amsmath}
\usepackage{amscd}
\usepackage{amssymb}
\usepackage{amsthm}
\usepackage{xspace}
\usepackage[all,tips]{xy}
\usepackage[dvips]{graphicx}
\usepackage{verbatim}
\usepackage{syntonly}
\usepackage{hyperref}

\theoremstyle{plain}
\newtheorem{thm}{Theorem}[section]
\newtheorem{cor}[thm]{Corollary}
\newtheorem{prop}[thm]{Proposition}
\newtheorem{lemma}[thm]{Lemma}

\theoremstyle{definition}

\newenvironment{pf}
{\begin{proof}} {\end{proof}}

\DeclareMathOperator{\SL}{SL} 
\DeclareMathOperator{\GL}{GL} 
 
 \DeclareMathOperator{\SO}{SO}

\DeclareMathOperator{\SU}{SU} 
\DeclareMathOperator{\Sp}{Sp} 
 
 \DeclareMathOperator{\Comm}{Comm}

 \DeclareMathOperator{\Mat}{M}

\DeclareMathOperator{\D}{D}\DeclareMathOperator{\Bo}{B}
\DeclareMathOperator{\Co}{C}\DeclareMathOperator{\A}{A}
\DeclareMathOperator{\Isom}{Isom}
\DeclareMathOperator{\SD}{SD}\DeclareMathOperator{\Fo}{F}

\newcommand{\vp}{\varphi}

\newcommand{\nid}{\noindent}

\newcommand{\iny}{\infty}

\newcommand{\co}{\ensuremath{\colon}}

\newcommand{\abs}[1]{\left\vert#1\right\vert}
\newcommand{\set}[1]{\left\{#1\right\}}

\newcommand{\pr}[1]{\left( #1 \right) }
\newcommand{\su}{\subset}

\newcommand{\lra}{\longrightarrow}

\newcommand{\B}[1]{\ensuremath{\mathbf{#1}}}

\newcommand{\Cal}[1]{\ensuremath{\mathcal{#1}}}

\newcommand{\Hy}{\ensuremath{\B{H}}}

\newcommand{\Q}{\ensuremath{\B{Q}}}
\newcommand{\R}{\ensuremath{\B{R}}}

\newcommand{\C}{\ensuremath{\B{C}}}
\newcommand{\F}{\ensuremath{\B{F}}}

\usepackage{fancyhdr}

\fancypagestyle{plain}

\begin{document}


\title{\textbf{Isospectral locally symmetric manifolds}}
\author{D. B. McReynolds}
\maketitle


\begin{abstract}
\nid In this article we construct closed, isospectral, non-isometric locally symmetric manifolds. We have three main results. First, we construct arbitrarily large sets of closed, isospectral, non-isometric manifolds. Second, we show the growth of size these sets of isospectral manifolds as a function of volume is super-polynomial. Finally, we construct pairs of infinite towers of finite covers of a closed manifold that are isospectral and non-isometric at each stage.
\end{abstract}

\nid keywords: \emph{isospectral tower, simple Lie group, Sunada's method, symmetric space.}\smallskip

\nid MSC code: 58J53, 11F06

\section{Introduction and main results}

\noindent For a closed Riemannian $n$--manifold $M$ with associated Laplacian $\Delta$, the eigenvalue spectrum of $\Delta$ acting on $L^2(M)$ is discrete with each eigenvalue occurring with finite multiplicity. We say two closed Riemannian $n$--manifolds $M_1$ and $M_2$ are \emph{isospectral} if $M_1$ and $M_2$ have equal eigenvalue spectra for their respective Laplacians. There is a long history on constructing isospectral non-isometric manifolds and for brevity we only touch on those results most pertinent to the present paper (see \cite{Gordon2000} for a survey). Sunada \cite{Sunada85} produced the first general (algebraic) method for constructing isospectral manifolds, and with this method, constructed many isospectral, non-isometric Riemann surfaces. Using \cite{Sunada85} and the Strong Approximation Theorem, Spatzier \cite{Spatzier89} showed every compact, irreducible, locally symmetric manifold admits a pair of isospectral non-isometric finite covers provided the isometry group of the universal cover has a complexification of sufficiently high rank. Using \cite{Sunada85}, Brooks, Gornet, and Gustafson \cite{BGG98} constructed arbitrarily large sets of isospectral non-isometric Riemann surfaces and studied the growth of the size of these sets as a function of volume (or equivalently, genera). \smallskip\smallskip

\noindent Our first result shows the existence of large sets of closed, isospectral, non-isometric manifolds modeled on the symmetric space $X$ associated to non-compact simple Lie group; occasionally we refer to such manifolds as $X$--manifolds

\begin{thm}\label{GeneralSequence}
Let $G$ be a non-compact simple Lie group with associated symmetric space $X$. Then for any $n$, there exist $n$ closed, isospectral, non-isometric manifolds with universal cover $X$.
\end{thm}

\noindent The following immediate corollary of Theorem \ref{GeneralSequence} was also previously unknown.

\begin{cor}\label{Iso}
Let $G$ be a non-compact simple Lie group with associated symmetric space $X$. Then there exist closed, isospectral, non-isometric manifolds with universal cover $X$.
\end{cor}

\noindent Theorem \ref{GeneralSequence} provides new examples for all $G$. Even for hyperbolic $n$--manifolds, which have garnered much attention (\cite{Bergeron00}, \cite{ChinburgFriedman99}, \cite{Reid92}, \cite{Vigneras80A}), Theorem \ref{GeneralSequence} provides many new examples. Moreover, for exceptional simple Lie groups, both Theorem \ref{GeneralSequence} and Corollary \ref{Iso} are new. \smallskip\smallskip

\noindent For a simple non-compact Lie group $G$ with associated symmetric space $X$ and any real positive $t$, let $\SD_X(t)$ denote the cardinality of the largest set of isometry classes of closed $X$--manifolds that are pairwise isospectral and have volume no greater than $t$. According to work of Pesce \cite{Pesce97}, $\SD_X(t)$ is finite for all $t$. Our next result provides nontrivial lower bounds for $\SD_X(t)$.

\begin{thm}\label{SpectralDiameter}
Let $G$ be a non-compact simple Lie group with associated symmetric space $X$. Then for every positive integer $r$, there exists a strictly increasing sequence $\set{t_j}$ such that
\[ \SD_X(t_j) \geq t_j^r. \]
\end{thm}

\noindent In most settings this function was not previously known to be unbounded. However, for Riemann surfaces, Brooks--Gornet--Gustafson \cite{BGG98} provided lower bounds on the order of $t^{\log t}$. Our next result generalizes this lower bound when $X$ is real hyperbolic $n$--space or complex hyperbolic $2$--space.

\begin{thm}\label{SpectralDiameterNonArithmetic}
Let $X$ be complex hyperbolic $2$--space or real hyperbolic $n$--space. Then there exists a constant $D$ and a strictly increasing sequence $\set{t_j}$ such that
\[ \SD_X(t_j) \geq t_j^{D\log(t_j)}. \]
\end{thm}

\noindent Before stating our next result, we require some additional terminology. Given two manifolds $M$ and $N$, we say the pair possesses an \emph{isospectral tower} if there exist two infinite towers of finite covers $\set{M_j}$ and $\set{N_j}$ of $M$ and $N$ such that the covers $M_j$ and $N_j$ are isospectral and non-isometric for all $j$. A manifold $M$ is said to possess an \emph{isospectral tower} if $(M,M)$ possesses an isospectral tower.

\begin{thm}\label{GeneralTower}
Let $G$ be a non-compact simple Lie group with associated symmetric space $X$. Then there exist closed manifolds with universal cover $X$ that possess isospectral towers.
\end{thm}

\noindent To our knowledge the only other tower constructions are given by Vig\'{e}rnas \cite{Vigneras80A} for certain closed manifolds modeled on $\SL(2,\R)$, $\SL(2,\C)$, and products of these groups, and by Lubotzky--Samuel--Vishne \cite{LSV05} who constructed isospectral towers for certain pairs of closed locally symmetric manifolds modeled on $\SL(n,\R)$ and $\SL(n,\C)$ for $n>2$. Even more remarkable is the fact that the towers constructed in \cite{LSV05} are for incommensurable manifolds $M,N$ (see \cite{CHLR}, \cite{PrasadRap}, and \cite{Reid92}), an unobtainable feature via Sunada's method. \smallskip\smallskip

\nid Finally it is well known (see \cite{Spatzier89}) that the above isospectral, non-isometric manifolds afford us the following corollary of Theorem \ref{GeneralSequence}.

\begin{cor}
Let $G$ be a non-compact simple Lie group. Then for every $n$ there exist $n$ measurably distinct properly ergodic actions of $G$ with equal, discrete spectra.
\end{cor}

\paragraph{Acknowledgments.}

I am grateful to Daniel Allcock, Misha Belolipetsky, Emmanuel Breuillard, Benson Farb, Skip Garibaldi, John Hammond, Inkang Kim, Alex Lubotzky, Jason Manning, Dave Morris, Alan Reid, Ralf Spatzier, and Matthew Stover for interesting and simulating conversations on the material of this article. In addition, the author was partially supported by an NSF postdoctoral fellowship and a Clay lift-off fellowship.

\paragraph{Notation.}

For a number field $k$, we denote the set of places or valuations by $V(k)$. The archimedean places will be denoted by $V_\iny(k)$ and the finite places by $V_f(k)$. We refer to the places in $V_\iny(k)$ as either real or complex. Given a place $\nu$ in $V(k)$, we denote the completion by $k_\nu$, the characteristic of the associated residue field by $p_\nu$, and the cardinality of the associated residue field by $q_\nu$. The ring of $k$--integers will be denoted by $\Cal{O}_k$ and the completion with respect to $\nu$ by $\Cal{O}_\nu$.\smallskip\smallskip

\nid We say a pair of subgroups $\Gamma_1,\Gamma_2$ in a group $G$ are \emph{commensurable} if $\Gamma_1\cap\Gamma_2$ is finite index in both $\Gamma_1$ and $\Gamma_2$. For a subgroup $\Gamma$ of a group $G$, we define the $G$--commensurator of $\Gamma$ to be
\[ \Comm_G(\Gamma) = \set{g \in G ~:~g^{-1}\Gamma g,\Gamma \text{ are commensurable}}. \]

\section{Sunada's method}

\nid We begin with a review of Sunada's method and Heisenberg groups over finite fields.

\paragraph{1. Sunada's theorem.}

For a pair of subgroup $H,K$ of a finite group $G$, we say $H$ and $K$ are \emph{almost conjugate} if for all $g$ in $G$, we have the equality
\[ \abs{H \cap [g]} = \abs{K \cap [g]}, \]
where $[g]$ denotes the $G$--conjugacy class of $g$. The following theorem of Sunada \cite{Sunada85} is our sole tool for producing isospectral covers; see \cite{LMNR} for a recent variant.

\begin{thm}[Sunada's theorem]\label{Sunada}
Let $M$ be a closed Riemannian $n$--manifold such that there exists a surjective homomorphism
\[ \vp\co \pi_1(M) \lra G. \]
If $H,K$ are almost conjugate subgroups of $G$ and $M_H,M_k$ are the finite covers of $M$ associated to the finite index subgroups $\vp^{-1}(H),\vp^{-1}(K)$, then $M_H,M_K$ are isospectral (length isospectral).
\end{thm}

\nid Even when $H,K$ are non-conjugate subgroups of $G$, the covers $M_H,M_K$ can still be isometric.

\paragraph{2. Heisenberg groups over finite rings.}

In order to employ Theorem \ref{Sunada}, we require large families of pairwise almost conjugate, non-conjugate subgroups. In \cite{BGG98}, such families were produced inside Heisenberg groups over finite rings. To this end, recall that for a commutative ring $R$ with identity, the \emph{Heisenberg group} $N_3(R)$ is the subgroup of $\SL(3,R)$ of upper triangular unipotent matrices. We call the subgroup
\[ H(R) = \set{\begin{pmatrix} 1 & x & 0 \\ 0 & 1 & 0 \\ 0 & 0 & 1 \end{pmatrix}} \]
the \emph{horizontal subgroup}.\smallskip\smallskip

\nid For our goal of producing large families of pairwise almost conjugate, non-conjugate subgroups, we work with finite fields $\F_q$ of degree $m$ over $\F_p$. Viewing $\F_q$ as an $m$--dimensional $\F_p$--vector space, we have two families of transformations. We have the family of $\F_p$--linear transformations $\Mat(m,\F_p)$ and the family of $\F_q$--linear transformations $\F_q$. Upon selecting an $\F_p$--basis for $\F_q$, we are afforded an injective homomorphism $\F_q \to \Mat(n,\F_p)$.\smallskip\smallskip

\nid Given a transformation $f$ in $\Mat(m,\F_p)$, we associate to $f$ the subgroup
\[ H_f = \set{\begin{pmatrix} 1 & x & f(x) \\ 0 & 1 & 0 \\ 0 & 0 & 1 \end{pmatrix}}. \]
We call $H_f$ the \emph{$f$--twisted horizontal subgroup}.

\begin{prop}\label{TwistedProposition}
Let $\F_q/\F_p$ be a degree $m$ extension and $f,g$ be elements of $\Mat(m,\F_p)$. Then $H_f,H_g$ are almost conjugate in $N_3(\F_q)$. In addition, $H_f,H_g$ are conjugate in $N_3(\F_q)$ if and only if $f-g$ is an element of $\F_q$.
\end{prop}

\nid Proposition \ref{TwistedProposition} produces $p^{m(m-1)}$ pairwise almost conjugate, non-conjugate subgroups of $N_3(\F_q)$. Viewing $N_3(\F_q)$ as a subgroup of the group $\GL(n,\F_q)$ for $n\geq 3$, it is conceivable that some of the subgroups above become conjugate. However, up to $\GL(n,\F_q)$--conjugation, we are left with at least $p^{m(m-1)-n^2}$ distinct conjugacy classes in $\GL(n,\F_q)$.

\paragraph{3. Broad plan.}

To produce isospectral manifolds we will use the subgroups $H_f$ of $N_3(\F_q)$. The utilization of these groups entails a few things. First, we must produce homomorphisms onto finite groups containing $N_3(\F_q)$. The existence of these homomorphisms follows from the Strong Approximation Theorem. In addition, we must ensure that we have enough distinct subgroups $H_f$ as already we have seen if $m\leq n$, we cannot guarantee the subgroups are non-conjugate. To achieve this, we will use certain lattices defined over number fields with desirable properties that via the Cebotarev Density Theorem suffice for ensuring that $m$ is sufficiently large. Last, we need to make sure the manifolds associated to these groups are not isometric. To control isometries between the various isospectral covers, we employ a recent result of Belolipetsky and Lubotzky \cite{BelolipetskyLubotzky05} in the arithmetic setting and Margulis' dichotomy \cite{Margulis91} in the non-arithmetic setting.\smallskip\smallskip

\nid Our approach was inspired by a beautiful paper of Belolipetsky and Lubotzky \cite{BelolipetskyLubotzky04} on isometry groups of hyperbolic $n$--manifolds (see \cite{McRSurvey} for a more in-depth discussion). In both this article and \cite{BelolipetskyLubotzky04}, the primary obstruction is unexpected isometries of finite covers. For Riemann surfaces, controlling these isometries is done via a dimensional argument in the moduli space of genus $g$ curves. Due to rigidity, one cannot hope to generalize this approach. The method taken in \cite{BelolipetskyLubotzky04} can be viewed as a discrete version of the dimensional argument used for Riemann surfaces; in \cite{BelolipetskyLubotzky04} the core argument is a counting argument. Our approach is similar. Namely, we will produce large families of isospectral manifolds and by a counting argument conclude that most are non-isometric.

\section{Heisenberg groups in non-compact simple groups}

\nid In order to use the almost conjugate subgroups from Proposition \ref{TwistedProposition}, we must produce homomorphisms onto finite groups that contain $N_3(\F_q)$. The desired homomorphisms follow from the Strong Approximation Theorem. That these finite groups contain $N_3(\F_q)$ requires some elementary results on non-compact simple Lie groups. Towards this latter goal, in this short section, we show that most non-compact simple Lie groups contain $N_3(\R)$.\smallskip\smallskip

\nid The classical non-compact simple Lie groups are the following families of groups:
\begin{align*}
\A_n:& \quad \SL(n-1,\R), \SL(n-1,\C), \SL((n-1)/2,\Hy),\SU(p,q), \quad p+q=n-1. \\
\Bo_n:& \quad \SO(2n+1,\C), \SO(p,q), \quad p+q = 2n+1. \\
\Co_n:& \quad \Sp(2n,\R), \Sp(2n,\C),\Sp(p,q), \quad p+q=n. \\
\D_n:& \quad \SO(2n,\C), \SO(n,\Hy), \SO(p,q), \quad p+q=2n.
\end{align*}

\begin{lemma}\label{HeisenbergSubgroupReals}
Let $G$ be a non-compact simple Lie group such that $G$ is not locally isomorphic to $\SO(n,1)$. Then $G$ contains a subgroup isomorphic to $N_3(\R)$.
\end{lemma}

\begin{pf}
It is well known that $N_3(\R)$ is a subgroup of $\SU(2,1)$, $\SL(3,\R)$, and $\Sp(4,\R)$. As any non-compact simple Lie group with real rank at least two contains either a copy of $\SL(3,\R)$ or $\Sp(4,\R)$ \cite[Proposition 1.6.2]{Margulis91}, we obtain $N_3(\R)$ subgroups for most non-compact simple Lie groups. The remaining groups $\SO(n,1),\SU(n,1),\Sp(n,1)$, and $\Fo_{4(-20)}$, aside from $\SO(n,1)$, contain a copy of $\SU(2,1)$. Note that $\SU(1,1)$ and $\Sp(1,1)$ are of type $\SO(n,1)$.
\end{pf}

\nid As all of the above inclusions are defined over $\Q$, these groups contain a $\Q$--defined copy of $N_3(\R)$.

\section{Lattices in non-compact simple Lie groups}

\nid In this section we construct families of lattices in non-compact simple Lie groups that will be used to produce isospectral manifolds whose universal cover is the symmetric associated to these Lie groups. The main feature we require is that these groups be defined over fields that contain a large prime order cyclic extensions of $\Q$. The existence of such forms follows from Borel--Harder \cite{BorelHarder}. However, for the reader's convenience, we give an explicit construction for classical groups.

\subsection{Algebraic prerequisites}

\nid Given an odd prime $\ell$, the cyclotomic field given by adjoining to $\Q$ an $\ell$--primitive root of unity will be denoted by $\Q(\zeta_\ell)$. The totally real subfield $\Q(\cos(2\pi/\ell))$ for which $\Q(\zeta_\ell)/\Q(\cos(2\pi/\ell))$ is a quadratic extension will also be of some use. By taking prime divisors of either $\ell-1$ or $(\ell-1)/2$, we can produce totally imaginary or totally real cyclic extensions of $\Q$ of arbitrarily large prime degree. \smallskip\smallskip

\nid We will also need number fields with a certain number of real and complex places. Let $F$ be a totally real cyclic extension of $\Q$ of prime degree $\ell$ and let $\lambda$ be a primitive generator for the extension. Set $\set{\lambda_1,\dots,\lambda_\ell}$ to be the set of Galois conjugates of $\lambda$ ordered by cardinality. Select a rational number $r$ between $\lambda_j$ and $\lambda_{j+1}$. For most $r$, the number field $F(\sqrt{\lambda-r})$  is a quadratic extension of $F$ with precisely $j$ complex places and $2(\ell - j)$ real places.\smallskip\smallskip

\nid For a number field $F$, an  $F$--quaternion algebra $D$, and a place $\nu$ in $V(F)$, we say that $D$ is \emph{ramified} at $\nu$ if $D \otimes_F F_\nu$ is a division algebra and \emph{split} if $D \otimes_F F_\nu$ is isomorphic to $\Mat(2,F_\nu)$. The Albert--Hasse--Brauer--Noether Theorem provides the existence of quaternion algebras with specified behavior at each place.\smallskip\smallskip

\nid Every $F$--quaternion algebra comes equipped with a pair of involution $\tau_c,\tau_r$ called quaternionic conjugation and reversion. In addition, if $D$ is an $E$--quaternion algebra and $E/F$ is a quadratic Galois extension, the non-trivial Galois involution $\tau$ on $E$ has an extension to $D$ that, in an abuse of notation, we denote by $\tau$.\smallskip\smallskip

\nid Given an involution $\tau$ on a cyclic, central, simple $F$--algebra $D$, there is an extension to $\Mat(n,D)$, the $F$--algebra of $n$ by $n$ matrices with coefficients in $D$. Specifically, the extension is defined to be the composition of matrix transposition (i.e., taking the transpose) and applying $\tau$ coefficient-wise. We say $h$ in $\Mat(n,D)$ is $\tau$--hermitian if $\tau(h) = h$. For instance, when $D$ is an $F$--quaternion algebra and $\tau$ is either $\tau_c$ or $\tau_r$, we have a notation of $\tau_c$ and $\tau_r$--hermitian forms.

\subsection{Families of lattices in non-compact simple Lie groups}

\nid We refer the reader to \cite{Morris} for a much more complete description of these lattices (see also \cite{Tits66}).

\paragraph{1. Lattices of type $\A_n$.}

We begin with lattices in the groups of type $\A_n$.

\begin{itemize}

\item $\SU(p,q)$. Let $F/\Q$ be a totally real cyclic extension of prime degree $\ell$. Let $E/F$ quadratic extension such that $E$ is totally imaginary with non-trivial Galois involution $\tau$. Let $h$ be a $\tau$--hermitian form  on $E^n$ such that $h$ has signature $(p,q)$ at one real place and signature $(n,0)$ otherwise. By Borel and Harish-Chandra \cite{BorelChandra62}, the group $\SU(h,\Cal{O}_E)$ is a cocompact lattice in $\SU(p,q)$.

\item $\SL(n,\R)$. Let $F/\Q$ be a totally real cyclic extension of prime degree $\ell$. Let $E/F$ be quadratic extension with $\ell-1$ complex places, 2 real places, and with non-trivial Galois involution $\tau$. Let $h$ be a $\tau$--hermitian form on $E^n$ such that $h$ has signature $(n,0)$ at every real place of $F$ where the signature is well defined. By Borel and Harish-Chandra \cite{BorelChandra62}, the group $\SU(h,\Cal{O}_E)$ is a cocompact lattice in $\SL(n,\R)$.

\item $\SL(n,\C)$. Let $F'/\Q$ be a totally real cyclic extension of prime degree $\ell$. Let $F/F'$ be a quadratic extension such that $F$ has exactly one complex place. Let $E/F$ a quadratic extension such that $E$ is totally imaginary with non-trivial Galois involution $\tau$. Let $h$ a $\tau$--hermitian form on $E^n$ such that $h$ has signature $(n,0)$ at every real place of $F$. By Borel and Harish-Chandra \cite{BorelChandra62}, the group $\SU(h,\Cal{O}_E)$ is a cocompact lattice in $\SL(n,\C)$.

\item $\SL(n,\Hy)$. Let $F/\Q$ be a totally real cyclic extension of prime degree $\ell$. Let $E/F$ be quadratic extension with $\ell-1$ complex places, 2 real places, and with non-trivial Galois involution $\tau$. Let $D$ be an $E$--quaternion algebra such that $D$ is ramified at the real places of $E$. Let $h$ to be a $\tau$--hermitian form on $D^n$ such that $h$ has signature $(2n,0)$ at every real place of $F$ where the signature is well defined. By Borel and Harish-Chandra \cite{BorelChandra62}, the group  $\SU(h,\Cal{O}_D)$ is a cocompact lattice in $\SL(n,\Hy)$, where $\Cal{O}_D$ is a maximal order in $D$.

\end{itemize}

\paragraph{2. Lattices of type $\Bo_n$.}

In this subsection, we treat groups of type $\Bo_n$.

\begin{itemize}

\item $\SO(p,q)$. Let $F/\Q$ be a totally real cyclic extension of prime degree $\ell$. Let $b$ be a bilinear form on $F^{2n+1}$ such that $b$ has signature $(p,q)$ at one real place and signature $(2n+1,0)$ otherwise. By Borel and Harish-Chandra \cite{BorelChandra62}, the group $\SO(B,\Cal{O}_F)$ is a cocompact lattice in $\SO(p,q)$.

\item $\SO(2n+1,\C)$. Let $F/\Q$ a totally real cyclic extension of prime degree $\ell$. Let $E/F$ be a quadratic extension such that $E$ has exactly $1$ complex places and $2\ell-2$ real places. Let $b$ be a bilinear form on $E^{2n+1}$ such that $b$ has signature $(2n+1,0)$ at every place of $F$ where the signature is well defined. By Borel and Harish-Chandra \cite{BorelChandra62}, the group $\SO(b,\Cal{O}_E)$ is a cocompact lattice in $\SO(2n+1,\C)$.

\end{itemize}

\paragraph{3. Lattices of type $\Co_n$.}
In this subsection we treat the groups of type $\Co_n$.

\begin{itemize}

\item $\Sp(p,q)$. Let $F/\Q$ be a totally real cyclic extension of prime degree $\ell$. Let $D$ an $F$--quaternion algebra such that  $D$ is ramified at all the real places of $F$. Let $h$ be $\tau_c$--hermitian form on $D^n$ such that $h$ has signature $(p,q)$ at one place of $F$ and signature $(n,0)$ otherwise. By Borel and Harish-Chandra \cite{BorelChandra62}, the group $\SU(h,\Cal{O}_D)$ is a cocompact lattice in $\Sp(p,q)$, where $\Cal{O}_D$ is a maximal order in $D$.

\item $\Sp(2n,\R)$. Let $F/\Q$ be a totally real cyclic extension of prime degree $\ell$. Let $D$ an $F$--quaternion algebra such that $D$ is split at one place of $F$ and ramified otherwise. Let $h$ be a $\tau_c$--hermitian form on $D^n$ such that $h$ has signature $(n,0)$ at every place where $D$ is ramified. By Borel and Harish-Chandra \cite{BorelChandra62}, the group $\SU(h,\Cal{O}_D)$ is a cocompact lattice in $\Sp(2n,\R)$, where $\Cal{O}_D$ is a maximal order in $D$.

\item $\Sp(2n,\C)$. Let $F/\Q$ be a totally real cyclic extension of prime degree $\ell$. Let $E/F$ be quadratic extension such that $E$ has 1 complex place and $2\ell-2$ real places. Let $D$ be an $E$--quaternion algebra such that $D$ is ramified at every real place of $E$. Let $h$ be a $\tau_c$--hermitian form on $D^n$ such that $h$ has signature $(n,0)$ at every place where $D$ is ramified. By Borel and Harish-Chandra \cite{BorelChandra62}, the group $\SU(h,\Cal{O}_D)$ is a cocompact lattice in $\Sp(2n,\C)$, where $\Cal{O}_D$ is a maximal order in $D$.

\end{itemize}

\paragraph{4. Lattices of type $\D_n$.}

For the construction of lattices in $\D_n$, we can take the construction used for groups of type $\Bo_n$. We simply change the dimension from $2n+1$ to $2n$. The one case that remains is $\SO(n,\Hy)$.

\begin{itemize}

\item $\SO(n,\Hy)$. Let $F/\Q$ be a totally real cyclic extension of prime degree $\ell$. Let $D$ an $F$--quaternion algebra such that $D$ is ramified at exactly one real place. Let $h$ be a $\tau_r$--hermitian form on $D^n$ such that $h$ has signature $(2n,0)$ at every place where $D$ splits. By Borel and Harish-Chandra \cite{BorelChandra62}, the group $\SU(h,\Cal{O}_D)$ is a cocompact lattice in $\SO(n,\Hy)$, where $\Cal{O}_D$ is a maximal order in $D$.

\end{itemize}

\paragraph{5. Lattices of exceptional type.}

For lattices in exceptional non-compact simple Lie groups, by Borel and Harder \cite{BorelHarder}, there exists a $F$--form $\B{G}$ where $F_0/\Q$ is a cyclic extension of arbitrarily large prime degree $\ell$ with $F_0 \su F$. In addition, the group $\B{G}(\Cal{O}_F)$ is a cocompact lattice in $\B{G}(\R)$.

\paragraph{6. Lattices in the group $\SO(n,1)$.}

For lattices in $\SO(n,1)$, we can use the non-arithmetic manifolds constructed in \cite{GPS88}. According to \cite{Lubotzky96}, these manifolds also have fundamental groups with the property that there exists a finite index subgroup that admits a surjective homomorphisms onto a non-abelian free group, i.e. they have large fundamental groups. The group $\SU(2,1)$ also possesses non-arithmetic, cocompact, large lattices (see \cite{Livne81}). Alternatively, we can use the construction above for the groups $\SO(p,q)$.

\paragraph{7. Lattices in simply connected, connected groups.}

For each lattice $\Gamma$ constructed above, we have an associated algebraic group $\B{H}$ defined over a number field $F$. In turn, we obtain an associated simply connected, connected, absolutely simple algebraic group $\B{G}$ and an associated principal arithmetic group commensurable with $\Gamma$. In order to employ the Strong Approximation Theorem, we work with $\B{G}$ and the associated principal arithmetic lattice. By construction, $\B{G}$ is defined over a number field $F$ that contains a totally real cyclic extension of prime degree $\ell$, where $\ell$ can be taken to be arbitrarily large.

\paragraph{8. Lattices in isometry groups.}

The projection from $\B{G}$ to the connected component of the adjoint form $\B{G}'$ of $\B{G}$ yields a cocompact lattice in $\Isom(X)$, where $X$ is the associated symmetric space. It may be the case that a pair of lattices $\Gamma_1,\Gamma_2$ in $\B{G}$ are not conjugate in $\B{G}$ but have projections that are conjugate in $\Isom(X)$. For our purposes here, note that there exists a constant $C_X$ such that the projection of a lattice $\Gamma_0$ to $\Isom(X)$ can become conjugate to at most $C_X$ new subgroups.

\section{Proof of Theorem \ref{GeneralSequence}}

In this section we prove Theorem \ref{GeneralSequence}.

\subsection{Groups not locally isomorphic to $\SO(n,1)$}

\nid We start first with arithmetic lattices non-compact simple Lie groups not locally isomorphic to $\SO(n,1)$. We will treat the lattices $\SO(n,1)$ separately.

\paragraph{1. Maps onto finite groups.}

Let $\Gamma$ be a principal arithmetic lattice arising from one of the above arithmetic constructions and let $\B{G}$ be the associated connected, simply connected, absolutely simple $F$--algebraic group. For each place $\nu$ in $V_f(F)$, we have the homomorphism
\[ r_\nu\co \B{G}(\Cal{O}_\nu) \lra \B{G}(\Cal{O}_\nu/\pi\Cal{O}_\nu), \]
where $\pi\Cal{O}_\nu$ is the uniformizer ideal.  According to the Strong Approximation Theorem (see \cite{Nori87}, \cite{Pink00}, \cite{Weisfeiler84}), there exists a finite set of places $S$ in $V_f(F)$ such that for for all $\nu$ in $V_f(F) \smallsetminus S$, the homomorphism $r_\nu$ restricted to $\Gamma$ is surjective. We also have
\[ \B{G}(\Cal{O}_\nu/\pi\Cal{O}_\nu) \cong \B{G}(\F_{q_\nu}), \]
where $\F_{q_\nu}$ is the residue field.

\paragraph{2. Heisenberg groups in completions.}

Using Lemma \ref{HeisenbergSubgroupReals}, we now argue that the groups $\B{G}(\F_{q_\nu})$ contain the finite group $N_3(\F_{q_\nu})$ for almost all $\nu$. For all but finitely many place $\nu$ in $V(F)$, the group $\B{G}(F_\nu)$ is quasi-split \cite[Theorem 6.7]{PlatonovRapinchuk94} and thus possesses a $F_\nu$--defined Borel subgroup. In addition, the injections of $N_3(\R)$ given by Lemma \ref{HeisenbergSubgroupReals} are all defined over $\Q$. This pair of facts in tandem with \cite[Proposition 3.21]{PlatonovRapinchuk94} show that for all but finitely many places $\nu$, we have an injection $N_3(F_\nu)$ into $\B{G}(F_\nu)$. Hence, for all but finitely many places, $\B{G}(\F_{q_\nu})$ contains $N_3(\F_{q_\nu})$ as a subgroup.

\paragraph{3. Controlling the finite fields.}

In order to use the subgroup $N_3(\F_{q_\nu})$ to produce almost conjugate, non-conjugate pairs in $\B{G}(\F_{q_\nu})$, we need $q_\nu$ to be sufficiently large. As the lattices constructed above are all defined over a field $F$ with a cyclic subfield of prime degree $\ell$, the Cebotarev Density Theorem provides us with a positive density subset of places $V_\ell$ of $F$ such that for all $\nu$ in $V_\ell$, we have
\[ q_\nu \geq p_\nu^\ell. \]
We take a positive density subset $V_1$ of $V_\ell$ such that
\[ [\F_{q_\nu}:\F_{p_\nu}] = \ell_0 \]
is constant and  $\ell_0\geq \ell$. We order the places $V_1 = \set{\nu_j}$ according to the characteristic of the associated residue field. By removing any redundant places, we may assume each characteristic $p_\nu$ that is realized is uniquely realized.

\paragraph{4. Subgroups via pullbacks.}

For each $\nu$ in $V_1$, let $\set{H_f}$ denote the $\B{G}(\F_{q_\nu})$--conjugacy classes of twisted horizontal subgroups and $\set{\Gamma_f}$, the set of finite index subgroups of $\Gamma$ given by $r_\nu^{-1}(H_f)$, where we now define $r_\nu$ to have domain $\Gamma$. By Proposition \ref{TwistedProposition}, there are at least
\begin{equation}\label{TwistedCount1}
p_\nu^{\ell_0(\ell_0-1) - \ell_0\dim \B{G}}
\end{equation}
distinct $\Gamma$--conjugacy classes of subgroups in the set $\set{\Gamma_f}$. Note that for (\ref{TwistedCount1}), we require the well known fact $\abs{\B{G}(\F_{q_\nu})} \sim q_\nu^{\dim \B{G}}$ (see \cite[p. 131]{Steinberg}).

\paragraph{4. Isometries of distinct covers.}

The following result of Belolipetsky and Lubotzky \cite[Corollary 5.3]{BelolipetskyLubotzky05} is the main tool we use to control the $\Comm_{\B{G}}$--conjugacy classes of the subgroups $\Gamma_f$ constructed above.

\begin{prop}\label{BLIsometryControl}
Let $\Gamma$ be a principal arithmetic lattice, $\Gamma_0$ an index $n$ congruence subgroup containing the principal congruence subgroup $\ker r_\nu$. Then there exists a constant $x$ (depending on $\Gamma$) and a constant $C$ (depending on $\B{G}$ only as Lie group) such that the number of index $n$ subgroups of $\Gamma$ that are $\Comm_{\B{G}}(\Gamma)$--conjugate to $\Gamma_0$ is at most
\[ nx^C\cdot \abs{\Gamma/\ker r_\nu}. \]
\end{prop}

\nid Applying Proposition \ref{BLIsometryControl} to the groups $\Gamma_f$, we see that there are at most
\[ x^Cq_\nu^{2\dim \B{G}} \]
subgroups $\Gamma_g$ that are $\Comm_{\B{G}}(\Gamma)$--conjugate to $\Gamma_f$. Taking $q_\nu \gg x$, we see that
\[ x^cq_\nu^{2\dim \B{G}} < q_\nu^{2\dim \B{G} + c}. \]
Taking into account that we already reduced by the $\Gamma_0$--conjugation, we see that in total, we have at least
\[ p_\nu^{\ell_0(\ell_0-1) - \ell_0(3\dim \B{G} + c)} \]
subgroups $\Gamma_f$ that have distinct $\Comm_{\B{G}}(\Gamma)$--conjugacy classes. Taking the degree $\ell$ of the cyclic subfield $F_0$ in $F$, the field of definition for $\B{G}$, so that
\[ \ell(\ell-1) - \ell(3\dim \B{G} + c) > r, \]
and recalling that $\ell \leq \ell_0$, we produce at least $p_\nu^r$ distinct $\Comm_{\B{G}}$--conjugacy classes of subgroups $\Gamma_f$.\smallskip\smallskip

\paragraph{5. Proof of Theorem \ref{GeneralSequence}.}

Passing to $\Isom(X)$, where $X$ is the associated symmetric space, there is a fixed constant $C_X$ that depends only on $X$ such that the image of $\Gamma_f$ is conjugate to at most $C_X$ additional images of the subgroups $\Gamma_g$ in $\Isom(X)$. In particular, taking $\ell$ and $p_\nu$ such that
\[ C_X^{-1}p_\nu^{\ell(\ell-1) - \ell(3\dim \B{G} + c)} > n, \]
by Theorem \ref{Sunada}, we produce at least $n$ compact orbifolds that are isospectral and non-isometric. A few words are required on why the associated orbifolds are isospectral. For each place $\nu$ in $V_1$, we have a surjective homomorphism
\[ r_\nu\co \Gamma \lra \B{G}(\F_{q_\nu}). \]
Projecting to $\Isom(X)$, the image of $\ker r_\nu$ induces a map onto a finite group $Q_\nu$ and a homomorphism $\B{G}(\F_{q_\nu})$ to $Q_\nu$. The image of the twisted horizontal subgroups $H_f$ in $Q_\nu$ are still almost conjugate; this fact is true in general. Thus, by Theorem \ref{Sunada}, the orbifolds associated to the images of the subgroups $\Gamma_f$ are isospectral. By Mostow Rigidity, these orbifolds are pairwise non-isometric. Since the group $N_3(\F_q)$ is unipotent, it follows from the proof of \cite[Proposition 3.8]{McRCusp3} that the groups $\Gamma_f$ are torsion free for all but finitely many places $\nu$ in $V_1$. Thus, the orbifolds $M_f$ are in fact manifolds, and hence Theorem \ref{GeneralSequence} holds. \qed

\subsection{Lattices in $\SO(n,1)$}

\nid For lattices in $\SO(n,1)$, so long as $n+1 \geq 5$, the associated groups $\B{G}(\F_q)$ also contain $N_3(\F_q)$ as a subgroup. In particular, when $n+1 \geq 5$, we can apply the methods above for lattices in $\SO(n,1)$. Alternatively, since there exists a non-arithmetic, large lattice $\Gamma$ in $\SO(n,1)$ (see \cite{Lubotzky96}) and the groups $\SL(3,\F_q)$ are two generator, there exists a finite index subgroup $\Gamma_0$ of $\Gamma$ and surjective homomorphisms
\[ r_\ell\co \Gamma_0 \lra \SL(3,\F_{p^\ell}) \]
for all $\ell$. Taking the pullbacks of the subgroups $H_f$ as before, we produce at least $p^{\ell(\ell-1) - 8\ell}$ subgroups $\Gamma_f$. As before, we can make
\[ \ell(\ell-1)- 8\ell \]
arbitrarily large.\smallskip\smallskip

\nid To control the non-arithmetic isospectral manifolds, we utilize the deep dichotomy established by Margulis (see \cite[Theorem 1, p. 2]{Margulis91}).

\begin{thm}[Margulis' dichotomy]\label{MargulusIsometryControl}
If $\Gamma_0$ is non-arithmetic, then
\[ [\Comm_{\Isom(X)}(\Gamma_0):\Gamma_0] < \iny. \]
In particular, there are at most $[\Comm_{\Isom(X)}(\Gamma_0):\Gamma_0]$ additional subgroups $\Gamma_g$ that are $\Comm_{\Isom(X)}(\Gamma_0)$--conjugate to $\Gamma_f$.
\end{thm}

\nid If we take $\ell$ such that
\[ \frac{p^{\ell(\ell-1)-8\ell}}{[\Comm_{\Isom(X)}(\Gamma_0):\Gamma_0]} > n, \]
we produce at least $n$ subgroups $\Gamma_f$ that are distinct up to $\Isom(X)$--conjugacy. By Mostow Rigidity and Theorem \ref{Sunada}, we obtain at least $n$ isospectral, non-isometric closed orbifolds. To ensure that we have produced manifolds opposed to orbifolds, by Selberg's Lemma, we select $\Gamma_0$ to be torsion free.

\section{Isospectral growth}

\nid In this section, we prove Theorem \ref{SpectralDiameter} and Theorem \ref{SpectralDiameterNonArithmetic}. We begin with Theorem \ref{SpectralDiameter}, the easier of the pair.

\paragraph{Proof of Theorem \ref{SpectralDiameter}.}

To produce enough isospectral, non-isometric manifolds modeled on $X$, we utilize our freedom in choosing the field $F$ used in the production of the arithmetic lattices. As the volume of the manifolds produced above is no greater than $C_1p_\nu^{\ell_0\dim \B{G}}$ for a constant $C_1$ that depends only on $\Gamma$, by taking $\ell$ sufficiently large, we can ensure that $\SD_X(t_j)$ is larger than $t_j^r$ for any $r$ and a sequence of $t_j$. Specifically, we have at least
\[ C_X^{-1}p_\nu^{\ell_0(\ell_0-1) - \ell_0(3\dim \B{G} + c)} \]
distinct isospectral, non-isometric covers with volume no more than $C_1p_\nu^{\ell\dim \B{G}}$. Selecting $\ell$ such that
\begin{equation}\label{BigEll0}
\ell(\ell-1) - \ell(3\dim\B{G} + c) - r\ell\dim \B{G} - r > 0,
\end{equation}
we produce the desired growth for Theorem \ref{SpectralDiameter}. Indeed, for large $p_\nu$, from (\ref{BigEll0}) we see that
\begin{align*}
C_1^rp_\nu^{r\ell_0 \dim \B{G}} &< p_\nu^{r+ r\ell_0\dim \B{G}} \\
&< p_\nu^{\ell_0(\ell_0-1) \ell_0(3\dim\B{G} + c)},
\end{align*}
as needed. \qed

\paragraph{Proof of Theorem \ref{SpectralDiameterNonArithmetic}.}

For Theorem \ref{SpectralDiameterNonArithmetic}, the argument is slightly more involved. To prove this result, we take a finite index subgroup $\Gamma_0$ of $\Gamma$ that admits a surjective homomorphism onto $\SL(3,\F_p[x])$. Our proof mimics the proof of the same result for Riemann surfaces given by Brooks--Gornet--Gustafson \cite{BGG98}.\smallskip\smallskip

\nid Let $\Gamma$ be a large cocompact lattice in the isometry group $\Isom(X)$ where $X$ is either $\Hy_\R^n$ or $\Hy_\C^2$. For every $p$, the groups $\SL(3,\F_p[x])$ are finitely generated and so for any $p$, there exists a finite index subgroup $\Gamma_0$ of $\Gamma$ and a surjective homomorphism
\[ \rho\co \Gamma_0 \lra \SL(3,\F_p[x]). \]
For each $j$, we have a surjective homomorphism
\[ r_j\co \SL(3,\F_p[x]) \lra \SL(3,\F_p[x]/x^j). \]
The group $\SL(3,\F_p[x]/x^j)$ contains the Heisenberg group $N_3(\F_p[x]/x^j)$. According to \cite{BGG98}, there are
\[ p^{j(j-1)/2} \]
$N_3(\F_p[x]/x^j)$--conjugacy classes of twisted horizontal subgroups $H_f$. Identifying those which are conjugate in $\SL(3,\F_p[x]/x^j)$, we obtain at least
\[ p^{j(j-1)/2 - 9j} \]
$\Gamma_0$--conjugacy classes of pullbacks of twisted horizontal subgroups. As $\Gamma$ is non-arithmetic, by Proposition \ref{MargulusIsometryControl}, there exists a constant $D_p$ such that there are at least
\[ p^{j(j-1)/2-9j - D_p} \]
$\Isom(X)$--conjugacy distinct pullbacks of twisted horizontal subgroups.\smallskip\smallskip

\nid For a fixed prime $p$ and each $j$, let $t_j$ denote the volume of a manifold $M_f$ associated to a twisted horizontal subgroup $H_f$ of $N_3(\F_p[x]/x^j)$. Then
\[ t_j = C[\Gamma_0:\Gamma_f] \leq Cp^{9j} \leq p^{9j+\delta}, \]
where $C$ is the volume of the orbifold associated to $\Gamma_0$ and $\delta$ is given by $p^\delta > C$. With this, we have
\begin{align*}
t_j^{\log(t_j)} &\leq \pr{p^{9j+\delta}}^{\log(p^{9j+\delta})} \\
&= p^{(9j+\delta)^2\log p}.
\end{align*}
From this calculation, it is clear that we can pick a constant $D$ such that
\[ D\log p(9j + \delta)^2 < j(j-1)/2-9j - D_p \]
for all $j$, since $p$, $\delta$, and $D_p$ are constant. In particular,
\[ t_j^{D\log(t_j)} < p^{j(j-1)/2 - 9j - D_p} \leq \SD_X(t_j) \]
as needed. \qed

\section{Isospectral towers}

\nid In this section, we prove Theorem \ref{GeneralTower}.

\paragraph{1. Some twisted variants.}

To produce towers of isospectral manifolds, we introduce a variant of the twisted horizontal subgroups of $N_3$ used above. To this end, we recall some notation from the proof of Theorem \ref{GeneralSequence}. To begin, we have a number field $F$ and a positive density set of places $V_1$ such that the associated residue field $\F_{q_\nu}$ has cardinality $p_\nu^{\ell_0}$ for all $\nu$ in $V_1$, where $\ell_0 \geq \ell$. We select $\ell$ so that
\[ \ell(\ell-1) - \ell(3\dim \B{G} + c) > r \]
holds for some $r>1$, where $c$ is a constant that depends on $\B{G}$ and $\Gamma$. We also order $V_1$ by the cardinality of the associated residue field of the place. In addition, removing redundances if necessary, we may assume each characteristic represented is uniquely represented. \smallskip\smallskip

\nid For each $j$, by the Strong Approximation Theorem, we have a surjective homomorphism
\[ r_j\co \Gamma \lra \prod_{i=1}^j \B{G}(\F_{q_{\nu_i}}) = Q_j. \]
Set
\[ L_j = \prod_{i=1}^j N_3(\F_{q_{\nu_i}}) < Q_j, \]
and note that the subgroup $L_j$ has at least
\[ \prod_{i=1}^j p_{\nu_i}^{\ell(\ell-1) - \ell(3\dim\B{G}-c)} \]
almost conjugate subgroups for which the pullbacks are distinct up to $\B{G}$--conjugation. In addition, we have a tower of homomorphisms
\[ \xymatrix{ Q_j \ar[rr]^{\psi_j} && Q_{j-1} \ar[rr]^{\psi_{j-1}} && \dots \ar[rr]^{\psi_3} && Q_2 \ar[rr]^{\psi_2} && Q_1} \]
given by projection.

\paragraph{2. Trees of manifolds.}

Viewing the subscript $j$ as a level, we see that the level $j$ twisted horizontal subgroups map via $\psi_j$ to level $j-1$ twisted horizontal subgroups. Forgetting the possible identifications arising from $\B{G}$--conjugation and conjugation in $\Gamma$, the number of level $k>j$ subgroups sitting over a fixed level $j$ subgroup is independent of the subgroup. In particular, the (counting) probability that a level $k>j$ subgroup sits over a level $j$ subgroup is independent of the level $k$ and of the level $j$ subgroup.\smallskip\smallskip

\nid The associated tree of covers we obtain from this family of groups, ignoring possible identifications due to isometries, is also uniform. We also refer to the level of the associated manifold via the identification.

\paragraph{3. Proof of Theorem \ref{GeneralTower}.}

To prove Theorem \ref{GeneralTower}, we need only make a few additional observations aside from the ones already made above. First, the construction used to prove Theorem \ref{GeneralSequence} in this setting yields at least a pair of level 1 manifolds that are isospectral and non-isometric. As the primes $p_{\nu_j}$ are increasing, the construction also shows that the probability that a given pair of level $k$ manifolds are isometric tends to zero as $k$ tends to infinity. In particular, if we start with a pair of level 1 non-isometric manifolds, there exists a level $k_2$ such that the probability a pair of level $k_2$ manifolds sitting over our pair of level 1 manifolds are isometric is less than 1. Indeed, given any pair of non-isometric manifolds of level $k$, there exists a level $k'$ and a pair of level $k'$ manifolds sitting over the level $k$ pair that are non-isometric. Consequently, we can iteratively build the desired tower needed for Theorem \ref{GeneralTower}.\smallskip\smallskip

\nid For the sake of completeness, we give some of the relevant calculations here. Namely, we prove the following claim, which suffices for an iterative construction.\smallskip\smallskip

\nid \textbf{Claim.} \emph{Given a pair of non-isometric level $j$ manifolds $M_1,M_2$, there exists a pair of level $j'$ non-isometric manifolds $N_1,N_2$ covering $M_1,M_2$.}\smallskip\smallskip

\nid Ignoring identifications as before, we know that there are
\[ \prod_{i=1}^j p_{\nu_i}^{\ell_0(\ell_0-1)} \]
level $j$ manifolds. For a fixed level $j$ manifold $M$, the number of level $k$ manifolds covering $M$ is
\[ \prod_{i=j+1}^k p_{\nu_i}^{\ell_0(\ell_0-1)}. \]
Finally, we know that for a fixed level $k$ manifold $N$, there are at most
\[ C_Xx^c\prod_{i=1}^k p_{\nu_i}^{3\ell_0 \dim \B{G}} \]
level $k$ manifolds isometric to $N$. In particular, to prove the claim, it suffices to show for sufficiently large $k$ that
\begin{equation}\label{TowerInequality}
C_Xx^c\prod_{i=1}^k p_{\nu_i}^{3\ell_0 \dim \B{G}} < \prod_{i=j+1}^k p_{\nu_i}^{\ell_0(\ell_0-1)}.
\end{equation}
Indeed, if (\ref{TowerInequality}) holds, then for a fixed level $k$ manifold $N_1$ covering $M_1$, the number of level $k$ manifolds cover $M_2$ that are isometric to $N_1$ would be less than the number of such manifolds. Consequently, there must exist the desired $N_2$ needed to validate the claim. To verify that (\ref{TowerInequality}) holds for sufficiently large $k$, we set
\[ C_j = \prod_{i=1}^j p_{\nu_i}^{3\ell_0\dim\B{G}}. \]
The constant $C_j$ depends only on $j$. In particular,
\[ C_Xx^c\prod_{i=1}^k p_{\nu_i}^{3\ell_0 \dim \B{G}} = C_XC_jx^c\prod_{i=j+1}^k p_{\nu_i}^{3\ell_0\dim\B{G}}. \]
As
\[ \ell(\ell-1) - 3\ell \dim\B{G} > r, \]
we see that
\[ \frac{\prod_{i=j+1}^k p_{\nu_i}^{\ell_0(\ell_0-1)}}{\prod_{i=j+1}^k p_{\nu_i}^{3\ell_0\dim\B{G}}} > \prod_{i=j+1}^k p_{\nu_i}^r. \]
If we select $k$ such that
\[ \prod_{i=j+1}^k p_{\nu_i}^r > C_jC_Xx^c, \]
then (\ref{TowerInequality}) holds.
\qed\smallskip\smallskip

\nid The construction of isospectral towers for non-arithmetic, large lattices is identical except we use a fixed homomorphism onto $\SL(3,\F_p[x])$. The tower of group homomorphisms
\[ \dots \lra \SL(3,\F_p[x]/x^j) \lra \SL(3,\F_p[x]/x^{j-1}) \lra \dots \lra \SL(3,\F_p[x]/x) \]
replaces the product construction in the arithmetic setting. The proof of the above claim in the non-arithmetic case is immediate via Theorem \ref{MargulusIsometryControl}.


\providecommand{\bysame}{\leavevmode\hbox to3em{\hrulefill}\thinspace}

\noindent
Department of Mathematics, University of Chicago\\
email: {\tt dmcreyn@math.uchicago.edu}


\end{document}